\documentclass[11pt]{amsart}
\usepackage{amsmath,amsthm,amssymb,url}
\usepackage[all]{xy}

\begin{document}

\newtheorem{theorem}{Theorem}[section]
\newtheorem{lemma}[theorem]{Lemma}
\newtheorem{proposition}[theorem]{Proposition}
\newtheorem{corollary}[theorem]{Corollary}
\newtheorem{bigtheorem}{Theorem}

\theoremstyle{definition}
\newtheorem{definition}[theorem]{Definition}
\newtheorem{example}[theorem]{Example}
\newtheorem{formula}[theorem]{Formula}
\newtheorem{nothing}[theorem]{}

\theoremstyle{remark}
\newtheorem{remark}[theorem]{Remark}

\renewcommand{\arraystretch}{1.2}

\newcommand{\de}{\partial}

\newcommand{\desude}[2]{{\dfrac{\de #1}{\de #2}}}

\newcommand{\mapor}[1]{{\stackrel{#1}{\longrightarrow}}}
\newcommand{\ormap}[1]{{\stackrel{#1}{\longleftarrow}}}
\newcommand{\mapver}[1]{\Big\downarrow\vcenter{\rlap{$\scriptstyle#1$}}}

\newcommand{\binfty}{\boldsymbol{\infty}}
\newcommand{\bi}{\boldsymbol{i}}
\newcommand{\bl}{\boldsymbol{l}}

\renewcommand{\bar}{\overline}
\renewcommand{\Hat}[1]{\widehat{#1}}

\newcommand{\sA}{\mathcal{A}}
\newcommand{\Oh}{\mathcal{O}}
\newcommand{\sF}{\mathcal{F}}
\newcommand{\sH}{\mathcal{H}}
\newcommand{\sL}{\mathcal{L}}
\newcommand{\sB}{\mathcal{B}}
\newcommand{\sY}{\mathcal{Y}}

\newcommand{\Q}{\mathbb{Q}}
\newcommand{\K}{\mathbb{K}}
\newcommand{\Proj}{\mathbb{P}}

\newcommand{\DER}{{{\mathcal D}er}}
\newcommand{\Id}{\operatorname{Id}}

\newcommand{\dee}{q_+}

\newcommand{\ad}{\operatorname{ad}}
\newcommand{\MC}{\operatorname{MC}}
\newcommand{\Def}{\operatorname{Def}}
\newcommand{\Hom}{\operatorname{Hom}}
\newcommand{\End}{\operatorname{End}}
\newcommand{\Image}{\operatorname{Im}}
\newcommand{\image}{\operatorname{Im}}
\newcommand{\Der}{\operatorname{Der}}
\newcommand{\Mor}{\operatorname{Mor}}
\newcommand{\Cone}{\operatorname{Cone}}
\newcommand{\Aut}{\operatorname{Aut}}
\newcommand{\coker}{\operatorname{Coker}}

\newcommand{\Grass}{\operatorname{Grass}}
\newcommand{\Flag}{\operatorname{Flag}}
\newcommand{\Hoch}{\operatorname{Hoch}}


\title[$L_{\infty}$ structures on mapping
cones]{$\mathbf{L}_{\boldsymbol\infty}$ structures on mapping
cones}

\date{April 3, 2007}

\author{Domenico Fiorenza}
\address{\newline Dipartimento di Matematica ``Guido
Castelnuovo'',\hfill\newline
Universit\`a di Roma ``La Sapienza'',\hfill\newline
P.le Aldo Moro 5,
I-00185 Roma Italy.\hfill\newline}
\email{fiorenza@mat.uniroma1.it}
\urladdr{www.mat.uniroma1.it/\~{}fiorenza/}
\author{Marco Manetti}
\email{manetti@mat.uniroma1.it}
\urladdr{www.mat.uniroma1.it/people/manetti/}

\begin{abstract}
We show that the mapping cone of a morphism of differential graded
Lie algebras $\chi\colon L\to M$ can be canonically endowed
with an $L_\infty$-algebra structure which at the same time
lifts the Lie algebra structure on $L$ and the usual
differential on the mapping cone. Moreover, this structure
is unique up to isomorphisms of $L_\infty$-algebras.
The associated deformation functor coincides with the one
introduced by the second author in
\cite{semireg}.
\end{abstract}

\subjclass{17B70, 13D10}
\keywords{Differential graded Lie algebras, symmetric coalgebras,
$L_{\infty}$-algebras, functors of Artin rings}

\maketitle

\section*{Introduction}

There are several cases where the tangent and obstruction
 spaces of a deformation theory are the cohomology
groups of the mapping cone of a morphism $\chi\colon L\to M$ of
differential graded Lie algebras. It is therefore natural to ask
if there exists a canonical differential graded Lie algebra
structure on the complex $(C_{\chi},\delta)$, where
\[C_\chi=\oplus C_\chi^i,\qquad C_\chi^i=
L^i\oplus M^{i-1},\qquad \delta(l,m)=(dl,\chi(l)-dm),\]%
such that the projection $C_{\chi}\to L$ is a morphism of
differential graded Lie algebras.\\
In general we cannot expect the existence of a Lie structure: in
fact the canonical bracket
\[
l_1\otimes l_2 \mapsto [l_1,l_2];\qquad m\otimes l\mapsto
\frac{1}{2}[m,\chi(l)];\qquad m_1\otimes m_2 \mapsto
0
\]
satisfies the Leibniz rule with respect to the differential
$\delta$ but  not the Jacobi identity. However, the Jacobi
identity for this bracket holds up to homotopy, and so we can look
for  the weaker request of a canonical $L_{\infty}$ structure on
$C_{\chi}$.\\
More precisely, let $\K$ be a fixed characteristic zero base
field, denote by $\mathbf{DG}$ the category of differential graded
vector spaces, by $\mathbf{DGLA}$ the category of differential
graded Lie algebras, by $\mathbf{L}_{\infty}$ the category of
$L_{\infty}$ algebras and by $\mathbf{DGLA}^{\mathbf 2}$ the
category of morphisms in $\mathbf{DGLA}$. The four functors
\[\begin{array}{cl}
\mathbf{DGLA}\to \mathbf{L}_{\infty}& \text{Quillen
construction},\\
\mathbf{L}_{\infty}\to \mathbf{DG}& \text{forgetting higher
brakets},\\
\mathbf{DGLA}^{\mathbf 2}\to \mathbf{DG}&
\{L\mapor{\chi}M\}\mapsto C_{\chi},\\
\mathbf{DGLA}\to \mathbf{DGLA}^{\mathbf 2}& L\mapsto\{L\to
0\},\end{array}
\]%

give a commutative diagram
\begin{equation*}
\xymatrix{\mathbf{DGLA}\ar[d]\ar[r]&\mathbf{L_\infty}\ar[d]\\
\mathbf{DGLA}^{\mathbf 2}\ar[r]&\mathbf{DG}}
\end{equation*}

Our first result is
\begin{bigtheorem}\label{thm.existence}
There exists a functor
$\widetilde{C}\colon\mathbf{DGLA}^{\mathbf 2}\to
\mathbf{L}_{\infty}$ making the  diagram
\[\xymatrix{
\mathbf{DGLA}\ar[d]\ar[r]&\mathbf{L_\infty}\ar[d]\\
\mathbf{DGLA}^{\mathbf
2}\ar[r]\ar[ur]^{\widetilde{C}}&\mathbf{DG}}\] commutative. If
$\sF\colon\mathbf{DGLA}^{\mathbf 2}\to \mathbf{L}_{\infty}$ has
the same properties, then for every morphism $\chi$ of
differential graded Lie algebras, the $L_{\infty}$-algebra
$\sF(\chi)$ is isomorphic to $\widetilde{C}(\chi)$.
\end{bigtheorem}

In the above theorem, the functor $\widetilde{C}$ is explicitly
described. The linear term of the $L_\infty$-algebra
$\widetilde{C}(\chi)$ is by construction the differential
$\delta$ on $C_\chi$, and the quadratic part turns out to
coincide with the naive bracket described at the beginning
of the introduction. An explicit expression for the higher
brackets is given in Theorem~\ref{thm.coefficienti}.

The second main result of this paper is to prove that the
deformation functor $\Def_{\widetilde{C}(\chi)}$ associated with
the $L_{\infty}$ algebra $\widetilde{C}(\chi)$ is isomorphic to
the functor $\Def_{\chi}$ defined in \cite{semireg}.

Given $\chi\colon L\to M$ it is defined a functor
$\Def_{\chi}\colon \mathbf{Art}\to \mathbf{Set}$, where
$\mathbf{Art}$ is the category of local Artinian $\K$-algebras
with residue field $\K$:
\[ \Def_{\chi}(A)=
\frac{\left\{(x,e^a)\in (L^1\otimes\mathfrak{m}_A)\times
\exp(M^0\otimes \mathfrak{m}_A)\mid
dx+\frac{1}{2}[x,x]=0,\;e^a\ast\chi(x)=0\right\}}
{\text{gauge equivalence}},\]%
where $\ast$ denotes the gauge action in $M$, and where
$(l_0,e^{m_0})$ is defined to be gauge equivalent to
$(l_1,e^{m_1})$ if  there exists $(a,b)\in (L^0\oplus
M^{-1})\otimes \mathfrak{m}_A$ such that
\[ l_1=e^a\ast l_0,\qquad e^{m_1}=e^{db}e^{m_0}e^{-\chi(a)}.\]

\begin{bigtheorem}\label{thm.isomorphismfunctors}
In the notation above, for every morphism of differential graded
Lie algebras $\chi\colon L\to M$ we have
\[ \Def_{\widetilde{C}(\chi)}\simeq\Def_{\chi}.\]
\end{bigtheorem}

The importance of Theorem~\ref{thm.isomorphismfunctors} relies on
the fact that it allows  to study the functors $\Def_{\chi}$, which are
often naturally identified with geometrically defined functors, 
using the whole machinery of
$L_\infty$-algebras. 
In particular this gives, under some finiteness assumption, 
the construction and the 
homotopy invariance of the Kuranishi 
map \cite{fuka,GoMil2,K},  as well as the local 
description of corresponding extended moduli spaces.\\

As a final remark, we observe that $\mathbf{DGLA^2}$ is a full
subcategory of the category  $\mathbf{DGLA^\Delta}$ of
cosimplicial differential graded Lie algebras and the
generalization of Theorem~\ref{thm.existence} to
$\mathbf{DGLA^\Delta}$ (with $C_{\chi}$ replaced by the total
complex) is essentially straightforward, using the ideas of this
paper and Whitney-Dupont
operators \cite{dupont,whitney,getzler}.\\
Naturally, it would be extremely interesting for applications to
deformation theory to prove the analogue of
Theorem~\ref{thm.isomorphismfunctors} for cosimplicial DGLAs: at
the moment we are unaware of simple descriptions of deformation
functors associated to cosimplicial DGLAs.

\medskip
\textbf{Acknowledgment.} Our thanks to Jim Stasheff for precious
comments on the  version v1 of this paper. The  version v3 of this
paper was written while the second author was at Mittag-Leffler
Institute in Stockholm, during a special year on Moduli Spaces:
the author is grateful for the support received and
for the warm hospitality.\\

\textbf{Keywords and general notation.}
We assume that the reader is familiar with the notion and main
properties
of differential graded Lie algebras and
$L_{\infty}$-algebras (we refer to
\cite{fuka,grassi1,K,LadaMarkl,LadaStas,defomanifolds}
as introduction of such structures); however
the basic definitions are recalled in this paper in order to fix
notation and terminology.\\
For the whole paper, $\mathbb{K}$ is a fixed field of
characteristic 0 and $\mathbf{Art}$ is the category of local
Artinian $\K$-algebras with residue field $\K$. For
$A\in\mathbf{Art}$ we denote by $\mathfrak{m}_A$ the maximal ideal
of $A$.

\bigskip

\section{Conventions on graded vector
spaces}\label{sec.graded.conv.}
In this paper we will work with $\mathbb{Z}$-graded vector spaces;
we write a graded vector space as $V=\oplus_{n\in {\mathbb Z}}V^n$,
and call $V^n$ the degree $n$ component of $V$; an element $v$ of
$V^n$ is called a degree $n$ homogeneous element of $V$. We say
that the graded vector space
$V$ is concentrated in degree
$k$ if
$V^i=\{0\}$ for $i\neq k$. Morphisms between graded vector spaces are
linear degree preserving maps, i.e. a map $\varphi\colon V\to W$
is a collection of linear maps $\varphi^n\colon V^n\to W^n$. The
shift functor is defined as $(V[k])^i:=V^{i+k}$. We say that a
linear map $\varphi\colon V\to W$ is a degree $k$ map if it is a
morphism $V\to W[k]$, i.e., if it is a collection of linear maps
$\varphi^n\colon V^n\to W^{n+k}$. The set of degree $k$
liner maps from $V$ to $W$ will be denoted $\Hom^k(V,W)$.
\par
Graded vector spaces are a symmetric
tensor category with $(V\otimes
W)^k=\oplus_{i+j=k}V^i\otimes W^j$ and $\sigma_{V,W}\colon V\otimes
W\to W\otimes V$ given by $\sigma(v\otimes w):=(-1)^{\deg(v)\cdot
\deg(w)}w\otimes v$ on homogeneous elements.
    We adopt the convention according to which
degrees are `shifted on the left'. By this we mean that we
have a natural identification, called the \emph{suspension}
isomorphism,
$V[1]\simeq
\mathbb{K}[1]\otimes V$ where $\mathbb{K}[1]$ denotes the graded vector
space consisting in the field $\mathbb{K}$ concentrated in
degree $-1$. Note that, with this convention the canonical
isomorphism $V\otimes \mathbb{K}[1]\simeq V[1]$ is
$v\otimes 1_{[1]}\mapsto (-1)^{\deg(v)}v_{[1]}$. More in general we have the
following \emph{decalage} isomorphism
\begin{align*}
V_1[1]\otimes\cdots \otimes V_n[1]&\xrightarrow{\sim}
(V_1\otimes\cdots \otimes V_n)[n]\\
{v_1}_{[1]}\otimes \cdots \otimes{v_n}_{[1]}&\mapsto
(-1)^{\sum_{i=1}^n(n-i)\cdot\deg{v_i}}(v_1\otimes \cdots \otimes
v_n)_{[n]}.
\end{align*}
Since graded vector spaces are a symmetric category, for any graded
vector space $V$ and any positive integer $n$ we have a canonical
representation of the symmetric group $S_n$ on $\otimes^n V$. The
space of coinvariants for this action is called the $n$-th
symmetric power of $V$ and is denoted by the symbol $\odot^n V$.
For instance
\[
V\odot V=V\otimes V/( v_1\otimes
v_2-(-1)^{\deg(v_1)\deg(v_2)}v_2\otimes v_1).
\]
Twisting the canonical representation of $S_n$ on $\otimes^n
V$ by the alternating character $\sigma\mapsto (-1)^\sigma$ and
taking the coinvariants one obtains the $n$-th antisymmetric (or
exterior) power of $V$, denoted by $\wedge^n V$. For instance
\[
V\wedge V=V\otimes V/( v_1\otimes
v_2+(-1)^{\deg(v_1)\deg(v_2)}v_2\otimes v_1).
\]
By naturality of the decalage isomorphism, we have a commutative
diagram
\[
\xymatrix{
      \bigotimes^n(V[1]) \ar[rr]^{{\rm decalage}}
\ar[d]_{\sigma[1]} &  &
\left(\bigotimes^nV\right)[n] \ar[d]^{(-1)^\sigma
\sigma[n]}\\
\bigotimes^n(V[1])  \ar[rr]^{{\rm decalage}} & &
\left(\bigotimes^nV\right)[n]\\
    }
\]
and so the decalage induces a canonical isomorphism
\[
\bigodot^n(V[1])\xrightarrow{\sim}\left(\bigwedge^n V\right)[n].
\]
As with ordinary vector spaces, one can identify $\odot^nV$ and
$\wedge^n V$ with suitable subspaces of $\otimes^nV$, called the
subspace of symmetric and antisymmetric tensors respectively.

\begin{remark}\label{rem.suspension}
Using the natural isomorphisms
\[ \Hom^i(V,W[l])\simeq \Hom^{i+l}(V,W)\]
and the decalage isomorphism, we obtain natural identifications
\[ \operatorname{dec}\colon\Hom^i\left(\bigotimes^k V,W\right)\xrightarrow{\sim}
\Hom^{i+k-l}\left(\bigotimes^k(V[1]),W[l]\right),\]
where
\[ \operatorname{dec}(f)(v_{1[1]}\otimes \cdots\otimes v_{k[1]})=
(-1)^{ki+\sum_{j=1}^k(k-j)\cdot\deg(v_j)}f(v_1\otimes \cdots\otimes v_{k})_{[l]}.\]
By the above considerations
\[ \operatorname{dec}\colon\Hom^i\left(\bigwedge^k V,V\right)\xrightarrow{\sim}
\Hom^{i+k-1}\left(\bigodot^k(V[1]),V[1]\right).\]

\end{remark}

\bigskip
\section{Differential graded Lie algebras and $L_\infty$-algebras}
\label{sec.linfty.intro}

A differential graded Lie algebra (DGLA for short) is a Lie algebra
in the category of graded vector spaces, endowed with a compatible degree 1 differential.
More explicitly, it is the
data $(V,d,[\,,\,])$, where $V$ is a graded vector space, the Lie
bracket
\[
[\,,\,]\colon V\wedge V\to V
\]
satisfies the graded Jacobi identity:
\[
[v_1,[v_2,v_3]]=[[v_1,v_2],v_3]+(-1)^{\deg(v_1)\deg(v_2)}[v_2,
[v_1,v_3]],
\]
and  $d\colon
V\to V$ is a degree 1 differential which is a degree 1
derivation of the Lie bracket, i.e.,
\[
d[v_1,v_2]=[dv_1,v_2]+(-1)^{\deg(v_1)}[v_1,dv_2].
\]
Morphisms of DGLAs are morphisms of graded vector spaces which are
compatible with the differential and the bracket, namely
\begin{align*}
\varphi(dv)&=d\varphi(v)\qquad \text{($\varphi$ is a morphism of
differential complexes)}\\
\varphi[v_1,v_2]&=[\varphi(v_1),\varphi(v_2)].
\end{align*}
Via the
decalage isomorphisms one can look at the Lie bracket of a DGLA $V$
as to a morphism
\[
q_2^{}\in\Hom^1(V[1]\odot V[1],V[1]),\qquad q_2(v_{[1]}\odot
w_{[1]})=(-1)^{\deg(v)}[v,w]_{[1]},
\]
Similarly, the suspended differential $q_1=d_{[1]}={\rm
id}_{\K[1]}\otimes d$ is a morphism of degree 1
\[
q_1^{}\colon V[1]\to V[1],\qquad q_1(v_{[1]})=-(dv)_{[1]}.
\]
Up to the canonical bijective linear map $V\to V[1]$, $v\mapsto v_{[1]}$,
the suspended differential $q_1$ and the bilinear operation
$q_2$ are written simply as
\[
q_1(v)=-dv; \qquad
q_2(v\odot w)=(-1)^{\deg_V(v)}[v,w],
\]
i.e., ``the suspended differential is the opposite
differential and $q_2$ is the twisted Lie
bracket''.

   Define morphisms $q_k^{}\in \Hom^1(
\odot^k(V[1]), V[1])$ by setting $q_k^{}\equiv 0$, for
$k\geq 3$. The map
\[
Q^1=\sum_{n\geq 1} q_n\colon \bigoplus_{n\geq 1}\bigodot^n
V[1]\to V[1]
\]
extends to a coderivation of degree 1
\[ Q\colon \bigoplus_{n\ge 1}\bigodot^n V[1]\to \left(\bigoplus_{n\ge
1}\bigodot^n
V[1]\right)\]
on the reduced symmetric coalgebra cogenerated by $V[1]$, by the
formula
\begin{nothing}\label{not.codifferential}
\[Q(v_1\odot\cdots\odot v_n)=\sum_{k=1}^n\sum_{\sigma\in S(k,n-k)}
\varepsilon(\sigma)q_k(v_{\sigma(1)}\odot\cdots\odot v_{\sigma(k)})
\odot v_{\sigma(k+1)}\odot\cdots\odot v_{\sigma(n)},\]
\end{nothing}
where
$S(k,n-k)$ is the set of unshuffles and
$\varepsilon(\sigma)=\pm 1$ is the  \emph{Koszul sign},
determined by the relation in
$\bigodot^n V[1]$
\[ v_{\sigma(1)}\odot\cdots\odot v_{\sigma(n)}=\varepsilon(\sigma)
v_{1}\odot\cdots\odot v_{n}.\]
The axioms of differential graded Lie algebra are then equivalent
to $Q$ being a codifferential, i.e., $QQ=0$. This
description of differential graded Lie algebras in terms of
the codifferential
$Q$ is called the Quillen construction \cite{Qui}. By
dropping the requirement that
$q_k^{}\equiv 0$ for
$k\geq 3$ one obtains the notion of $L_\infty$-algebra (or strong
homotopy Lie algebra), see e.g. \cite{LadaMarkl,LadaStas,K};
namely, an
$L_{\infty}$ structure on a graded vector space
$V$ is a sequence
of linear maps of degree 1
\[ q_k\colon \bigodot^k V[1]\to
V[1],\qquad k\ge 1,\]
such that the induced coderivation $Q$
on the reduced symmetric coalgebra cogenerated by $V[1]$, given  by
the Formula~\ref{not.codifferential} is a codifferential,
i.e., $QQ=0$. This condition  implies
$q_1q_1=0$ and therefore an $L_\infty$-algebra is in particular a
differential complex. Note that, by the above
discussion, every DGLA can be naturally seen as an
$L_\infty$-algebra; namely, a DGLA is an $L_\infty$-algebra with
vanishing higher multiplications $q_k^{}$, $k\geq 3$. Via the
decalage isomorphisms of Remark~\ref{rem.suspension},
the multiplications $q_k^{}$ of an
$L_\infty$-algebra $V$ can be seen as morphisms
\[
[\,,\dots,\,]_n^{}\in\Hom^{2-n}(\bigwedge^n V, V).
\]
The condition $QQ=0$ then translates into a sequence of
quadratic relations between the brackets
$[\,,\dots,\,]_n^{}$, the first of which are
$[[v]_1^{}]_1^{}=0$, i.e., $[\,]_1^{}$ is a degree 1
differential;
$[[v_1,v_2]_2^{}]_1^{}=[[v_1]_1^{},v_2]+(-1)^{\deg_V(v_1)}
[v_1,[v_2]_2^{}]^{}_2$, i.e., $[\,]_1^{}$ is a degree 1
derivation of the bracket $[\,,\,]^{}_2$;
\begin{align*}[v_1,[v_2,v_3]_2^{}]_2^{}&-[[v_1,v_2]_2^{},v_3]_2^{}-
(-1)^{\deg_V(v_1)\deg_V(v_2)}[v_2,
[v_1,v_3]_2^{}]_2^{}\\
&=[[v_1,v_2,v_3]_3^{}]_1^{}+[[v_1]_1^{},v_2,v_3]_3^{}+(-1)^{\deg_V(v_1)}
[v_1,[v_2]_1^{},v_3]_3^{}\\
&\qquad+(-1)^{\deg_V(v_1)+\deg_V(v_2)}
[v_1,v_2,[v_3]_1^{}]_3^{},
\end{align*}
i.e. the bracket $[\,,\,]^{}_2$ satisfies the graded Jacobi
identity up to the $[\,]_1^{}$-homotopy $[\,,\,,\,]_3^{}$, which
explains the name `homotopy Lie algebras'.
Note in particular that the $[\,]_1^{}$ cohomology of an
$L_\infty$-algebra carries a natural structure of graded
Lie algebra (or differential graded Lie algebra with trivial
differential).
\par A morphism
$f_\infty^{}$ between two $L_\infty$-algebras
$(V,q_1^{},q_2{},q_3^{},\dots)$ and
$(W,\hat{q}_1^{},\hat{q}_2{},\hat{q}_3^{},\dots)$ is a
sequence of linear maps of degree 0
\[ f_n\colon \bigodot^n V[1]\to
W[1],\qquad n\ge 1,\]
such that the morphism of coalgebras
\[ F\colon \bigoplus_{n\ge 1}\bigodot^n V[1]\to \bigoplus_{n\ge
1}\bigodot^n
W[1]\]
induced by $F^1=\sum_n f_n\colon\bigoplus_{n\ge 1}\bigodot^n
V[1]\to W[1]$
commutes with the codifferentials induced by the two $L_{\infty}$
structures on $V$ and
$W$ \cite{fuka,K,LadaMarkl,LadaStas,defomanifolds}.
An $L_{\infty}$-morphism $f^{}_\infty$ is called
\emph{linear}  (sometimes \emph{strict}) if $f_n=0$ for every
$n\ge 2$.  We note that a linear map $f_1\colon V[1]\to W[1]$
is a linear $L_{\infty}$-morphism if and only if
\[ \hat{q}_n(f_1(v_1)\odot\cdots\odot
f_1(v_n))=f_1(q_n(v_1\odot\cdots\odot v_n)),\qquad
\forall\; n\ge 1,\; v_1,\ldots,v_n\in V[1].\]
The category of $L_\infty$-algebras will be denoted by
$\mathbf{L}_\infty$ in this paper. Morphisms between DGLAs are
linear morphisms between the corresponding $L_\infty$-algebras, so
the category of differential graded Lie algebras is a
(non full) subcategory of $\mathbf{L}_\infty$.
\par
If $f_\infty$ is an $L_\infty$ morphism between
$(V,q_1^{},q_2{},q_3^{},\dots)$ and
$(W,\hat{q}_1^{},\hat{q}_2{},\hat{q}_3^{},\dots)$, then its
linear part
\[
f_1\colon V[1]\to W[1]
\]
satisfies the equation $f_1\circ q_1=\hat{q}_1\circ f_1$,
i.e., $f_1$ is a map of differential complexes
$(V[1],q_1)\to (W[1],\hat{q}_1)$. An
$L_\infty$-morphism $f^{}_\infty$ is called a quasiisomorphism of
$L_\infty$-algebras if its linear part $f_1$ is a quasiisomorphism
of differential complexes.

\bigskip
\section{The suspended mapping cone of $\chi\colon L\to M$.}
\label{sec.suspendedmap}

The suspended mapping cone of the DGLA morphism $\chi\colon L\to
M$ is the graded vector space
\[
C_\chi=\Cone(\chi)[-1],
\]
where $\Cone(\chi)=L[1]\oplus M$ is the mapping cone of
$\chi$. More explicitly,
\[ C_\chi=\mathop{\oplus}_iC_\chi^i,\qquad C_\chi^i=
L^i\oplus M^{i-1}.\]
The suspended mapping cone has  a natural differential
$\delta\in\Hom^1(C_{\chi},C_{\chi})$
given by
\[ \delta(l,m)=(dl,\chi(l)-dm),\qquad l\in L,m\in M.\]%

Denote $M[t,dt]=M\otimes \K[t,dt]$ and define, for every $a\in
\K$, the evaluation morphism
\[ e_a\colon M[t,dt]\to M,\qquad e_a(\sum m_it^i+n_it^idt)=\sum m_ia^i.\]
It is easy to prove that every morphism $e_a$ is a surjective
quasi-isomorphism of DGLA.\\
Consider the DGLA
\[ H_{\chi}=\{(l,m)\in L\times M[t,dt]\mid e_0(m)=0,\, e_1(m)=\chi(l)\}\]
and the morphism
\[ \imath\colon C_{\chi}\to H_{\chi},\qquad \imath(l,m)=(l,t\chi(l)+dt\cdot m).\]

\begin{proposition}\label{prop.esistequasiiso}
In the above notation, the morphism $\imath$ is an injective
quasi-isomorphism of complexes. For every functor $\sF\colon
\mathbf{DGLA}^{\mathbf 2}\to \mathbf{L}_{\infty}$
such that the diagram
\[\xymatrix{
\mathbf{DGLA}\ar[d]\ar[r]&\mathbf{L_\infty}\ar[d]\\
\mathbf{DGLA}^{\mathbf
2}\ar[r]\ar[ur]^{\sF}&\mathbf{DG}}\]
commutes, there exists an $L_{\infty}$-morphism
\[ \imath_{\infty}\colon \sF(\chi)\to H_{\chi}\]
with linear term $\imath_1=\imath$.
\end{proposition}

\begin{proof}
Defining
\[ P=\{(l,m)\in L\times M[t,dt]\mid  e_1(m)=\chi(l)\}.\]
We have a commutative diagram of morphisms of differential graded
Lie algebras
\[ \begin{array}{cccccc}
L&\xrightarrow{f}&P&\xleftarrow{}&H_{\chi}&\qquad f(l)
=(l,\chi(l))\\
\mapver{\chi}&&\mapver{\eta}&&\mapver{}&
\qquad \eta(l,m)=e_0(m)\\
M&\xrightarrow{\Id_M}&M&\xleftarrow{}&0&
\end{array}\]
and then 
two $L_{\infty}$-morphisms
$\sF(\chi)\rightarrow\sF(\eta)
\xleftarrow{h_{\infty}}H_{\chi}$ whose linear parts are
the two  injective
quasiisomorphisms
\[ C_{\chi}\rightarrow C_{\eta}\xleftarrow{h}H_{\chi},\qquad
h(l,m)=((l,m),0).\]%
A morphism of complexes $p\colon C_{\eta}\to H_{\chi}$ such that
$ph=\Id_{H_{\chi}}$ can be defined as
\[ p((l,m),n)=(l,m+(t-1)e_0(m)+dt\cdot n).\]
The composition of $p$ with the injective quasi-isomorphism
$C_{\chi}\to C_{\eta}$ gives the map $\imath$. 
By  general theory there exists a (non canonical) left inverse of
$h_{\infty}$ with linear term equal to $p$ and then the morphism
$\imath\colon C_{\chi}\to H_{\chi}$ can be lifted to an
$L_{\infty}$-quasiisomorphism.
\end{proof}

Denote by $\langle\;\rangle_1\in \Hom^1(C_{\chi}[1],C_{\chi}[1])$
and $q_1\in \Hom^1(H_{\chi}[1],H_{\chi}[1])$ the suspended
differentials, namely
\[ \langle (l,m)\rangle_1=(-dl,-\chi(l)+dm),\qquad l\in L,m\in M.\]%
\[ {q}_1(l,m)=(-dl,-dm).\]
Notice that $\imath$ induce naturally an injective
quasiisomorphism
\[\imath\colon C_{\chi}[1]\to H_{\chi}[1],\qquad
\imath(l,m)=(l,t\chi(l)+dt\cdot m).\]%

    The integral operator $\int_a^b\colon
\mathbb{K}[t,dt]\to \mathbb{K}$ extends  to a linear map of degree
$-1$
\[ \int_a^b\colon M[t,dt]\to M,\qquad
\int_a^b(\sum_i t^im_i+t^idt\cdot n_i)=
\sum_i\left(\int_a^bt^idt\right)n_i.\]

\begin{lemma}\label{lem.definitionK}
In the above notation, consider the linear maps
\[ \pi\in\Hom^0(
H_\chi[1],C_\chi[1]),\qquad K\in\Hom^{-1}(H_\chi[1],H_\chi[1])\]%
defined as
\[
\pi(l,m(t,dt))=\left(l,\int_0^1 m(t,dt)\right),\qquad
K(l,m)=\left(0,\int_0^t m-t\int_0^1m\right).
\]
Then $\pi$ is a morphism of complexes and
\[ \pi\imath=\operatorname{Id}_{C_\chi[1]},\qquad
\imath\,\pi=\operatorname{Id}_{H_\chi[1]}+Kq_1+{q_1}K.
\]
\end{lemma}

\begin{proof} The proof of the Lemma is a simple exercise;
we leave it to the reader.
\end{proof}

\bigskip
\section{Homotopy transfer of $L_{\infty}$ structures}
\label{sec.homotopytransfer}

A major result in the theory of $L_\infty$-algebras is the
following \emph{homotopical transfer of structure} theorem, which
we learnt from \cite{fuka,KonSoi}. For the reader's
convenience, we give a sketch of the proof. Note that the
version of the homotopy transfer we give here is slightly
more general than the ones we are aware of in the existing
litarature.

\begin{theorem}\label{thm.transfer}
Let $(V,q_1^{},q_2{},q_3^{},\dots)$ be an $L_\infty$-algebra and
$(C,\delta)$ be a differential complex. If there exist two
morphisms of differential complexes
\[
\imath\colon (C[1],\delta_{[1]}) \to (V[1],q_1) \qquad \text{and}
\qquad \pi\colon (V[1],q_1)\to (C[1],\delta_{[1]})
\]
such that the composition $\imath\pi$ is homotopic to the identity,
then there exist an $L_\infty$-algebra structure
$(C,\langle\,\rangle_1^{},\langle\,\rangle_2^{},\dots)$ on $C$
extending its differential complex structure and  an
$L_\infty$-morphism $\imath_\infty^{}$ extending $\imath$.
\end{theorem}

\begin{proof}
Let $K\in\Hom^{-1}(V[1],V[1])$ be an homotopy between $\imath\pi$
and $\operatorname{Id}_{V[1]}$, i.e.,
$q_1K+Kq_1=\imath\pi-\Id_{V[1]}$. Denote
\[\dee=\sum_{n\ge 2}q_n \colon \bigoplus_{n\geq
2}\bigodot^nV[1]\to V[1],\] so that  $Q^1=q_1+\dee$. Define a
morphism of  of graded coalgebras
\[
\imath^{}_\infty\colon \bigoplus_{n\geq 1}\bigodot^nC[1]\to
\bigoplus_{n\geq 1}\bigodot^nV[1]
\]
via the recursion
\[
\imath_\infty^1=\imath^1_1+K\dee\imath^{}_\infty
\]
and a degree 1 coderivation
\[
\hat{Q}\colon \bigoplus_{n\geq 1}\bigodot^nC[1]\to
\bigoplus_{n\geq 1}\bigodot^nC[1]
\]
by the formula
\[
\hat{Q}^1=\sum_{n\geq
1}\langle\,\rangle_n=\delta_{[1]}+\pi\dee\imath^{}_\infty.
\]
Then we have
\[ (Q\imath_\infty^{}-\imath_\infty^{}\hat{Q})^1=
 K\dee(Q\imath^{}_\infty-\imath^{}_\infty\hat{Q}).\]
Indeed,
\begin{equation*}
\begin{aligned}
(Q\imath^{}_\infty-\imath^{}_\infty\hat{Q})^1&=Q^1\imath^{}_\infty-
\imath_\infty^1\hat{Q}
=q_1\imath_\infty^1+\dee\imath^{}_\infty-\imath_\infty^1\hat{Q}\\
&=q_1\imath^1_1+q_1K\dee\imath^{}_\infty+\dee\imath^{}_\infty
-\imath^1_1\hat{Q}^1-K\dee\imath^{}_\infty\hat{Q}\\
&=q_1\imath+(\imath\pi-\Id_V
-Kq_1)\dee\imath^{}_\infty+\dee\imath^{}_\infty
-\imath\hat{Q}^1-K\dee\imath^{}_\infty\hat{Q}\\
&=q_1\imath+\imath\pi\dee\imath^{}_\infty -Kq_1\dee\imath^{}_\infty
-\imath\hat{Q}^1-K\dee\imath^{}_\infty\hat{Q}\\
&=q_1\imath+\imath\pi\dee\imath^{}_\infty-Kq_1\dee\imath^{}_\infty
-\imath\hat{Q}^1_1-\imath\pi\dee\imath^{}_\infty-K\dee\imath^{}_\infty\hat{Q}\\
&=(q_1\imath-\imath\delta_{[1]})-Kq_1\dee\imath^{}_\infty
-K\dee\imath^{}_\infty\hat{Q}\\
&=-Kq_1\dee\imath^{}_\infty-K\dee\imath^{}_\infty\hat{Q}.\\
\end{aligned}
\end{equation*}
Since $0=Q^1Q=q_1Q^1+\dee Q=q_1\dee+\dee Q$ we have $q_1\dee=-\dee
Q$ and therefore
\[
(Q\imath_\infty^{}-\imath_\infty^{}\hat{Q})^1
=K\dee(Q\imath_\infty{}-\imath_\infty^{}\hat{Q}).\]

The map
\[Q\imath_\infty^{}-\imath_\infty^{}\hat{Q}\colon \bigoplus_{n\geq
1}\bigodot^nC[1]\to \bigoplus_{n\geq 1}\bigodot^nV[1]\] is a
$\imath_\infty^{}$-derivation and then, in order to prove that
$Q\imath_\infty^{}-\imath_\infty^{}\hat{Q}=0$, it is sufficient to
show that $(Q\imath_\infty^{}-\imath_\infty^{}\hat{Q})^1=0$. We
shall prove by induction on $n$ that
$(Q\imath_\infty^{}-\imath_\infty^{}\hat{Q})^1$ vanishes on
$\bigodot^n C[1]$; for $n=0$ there is nothing to prove. Let us
assume $n>0$ and
$(Q\imath_\infty^{}-\imath_\infty^{}\hat{Q})^1(\bigodot^i C[1])=0$
for every $i<n$; then by coLeibniz rule, for every $w\in
\bigodot^n C[1]$ we have
$(Q\imath_\infty^{}-\imath_\infty^{}\hat{Q})(w)=
(Q\imath_\infty^{}-\imath_\infty^{}\hat{Q})^1(w)\in V[1]$ and
therefore
\[
(Q\imath_\infty^{}-\imath_\infty^{}\hat{Q})^1(w)=
K\dee(Q\imath_\infty^{}-\imath_\infty^{}\hat{Q})(w)=
K\dee(Q\imath_\infty^{}-\imath_\infty^{}\hat{Q})^1(w)=0.\] Therefore
\[
Q\imath_\infty^{}=\imath_\infty^{}\hat{Q}.
\]
We also have
\begin{align*}(\hat{Q}\hat{Q})^1&=\hat{Q}^1\hat{Q}=\hat{Q}^1_1\hat{Q}+
\pi\dee\imath_\infty^{}\hat{Q} =\delta_{[1]}\hat{Q}^1+
\pi\dee Q\imath_\infty^{}\\
&=\delta_{[1]}\pi\dee\imath_\infty^{}+\pi\dee Q\imath_\infty^{}=
\pi(q_1\dee+\dee Q)\imath_\infty^{}.
\end{align*}
We have already noticed that $q_1\dee=-\dee Q$ and then
$(\hat{Q}\hat{Q})^1=0$. Since $\hat{Q}$ is a coderivation, we find
\[
\hat{Q}\hat{Q}=0.
\]
\end{proof}

The recursive definition of $\imath^1_\infty$ can be explicitly
solved in terms of a summation over rooted trees \cite[Definition
6]{KoSo}; see also \cite{fuka,Schuhmacher}. Similarly, also the
operator $\hat{Q}^1$ can be written as a sum over rooted trees. We
sketch a proof of these facts following \cite{fiorenza:graphs}.
Let ${\mathcal T}_{K}$ be the groupoid whose objects are directed
rooted trees with internal vertices of valence at least two; trees
in ${\mathcal T}_{K}$ are  decorated as follows: each tail edge of
a tree in ${\mathcal T}_{K}$ is decorated by the operator
$\imath$, each internal edge is decorated by the operator $K$ and
also the root edge is decorated by the operator $K$; every
internal vertex $v$ carries the  operation $q_r$, where $r$ is the
number of edges having $v$ as endpoint. Isomorphisms between
objects in ${\mathcal T}_{K}$ are isomorphisms of the underlying
trees. Denote the set of isomorphism classes of objects of
${\mathcal T}_K$ by the symbol $T_K$. Similarly, let ${\mathcal
T}_{\pi}$ be the groupoid whose objects are directed rooted trees
with the same decoration as ${\mathcal T}_{K}$ except for the
root edge, which is decorated by the operator $\pi$ instead of
$K$. The set of isomorphism classes of objects of ${\mathcal
T}_\pi$ is denoted $T_\pi$.\\
Via the usual operadic rules, each decorated tree $\Gamma\in
{\mathcal T}_K$ with
$n$ tail vertices gives a linear map
\[ Z_\Gamma(\imath,\pi,K,q_i)\colon C[1]^{\odot n}\to V[1].\]
More precisely, let $\tilde{\Gamma}$ ba a planar
embedding of $\Gamma$; the standard orientation of the plane
induces a total ordering on tail vertices and then also a map
\[ Z_{\tilde{\Gamma}}(\imath,\pi,K,q_i)\colon C[1]^{\otimes n}\to
V[1]\] evalued according to the usual operadic rules. Then define
$Z_\Gamma(\imath,\pi,K,q_i)$ as the composition of
$Z_{\tilde{\Gamma}}(\imath,\pi,K,q_i)$ and the symmetrization map
\[C[1]^{\odot n}\to C[1]^{\otimes n},\quad v_1\odot\cdots\odot
v_n\mapsto \sum_{\sigma\in S_n}\epsilon(\sigma)
v_{\sigma(1)}\otimes\cdots\otimes v_{\sigma(n)}.\]%
It is straightforward to check that $Z_\Gamma(\imath,\pi,K,q_i)$
is well defined.\\
Similarly, each decorated tree in ${\mathcal T}_\pi$ gives rise to
a degree one multilinear operator on $C[1]$ with values in $C[1]$.
Here is an example:
\[
\Gamma=\qquad\begin{xy}
,(-10,6.66);(-6,4)*{\,\scriptstyle{\imath}\,}**\dir{-}
,(-10,-6.66);(-6,-4)*{\,\scriptstyle{\imath}\,}**\dir{-}
,(-6,4)*{\,\scriptstyle{\imath}\,};
(0,0)*{\,\,\scriptstyle{q_2}\,}**\dir{-}?>*\dir{>}
,(-6,-4)*{\,\scriptstyle{\imath}\,};(0,0)*{\,\,\scriptstyle{q_2}\,}**\dir{-}?>*\dir{>}
,(0,0)*{\,\,\scriptstyle{q_2}\,};
(8,0)*{\,\scriptstyle{\pi}\,}**\dir{-}
,(8,0)*{\,\scriptstyle{\pi}\,}; (16,0)**\dir{-}?>*\dir{>}
\end{xy},\]
\[Z_\Gamma(\imath,\pi,K,q_i)(a\odot b)=\pi q_2(i(a)\odot
i(b))+(-1)^{\bar{a}\;\bar{b}}\pi q_2(i(b)\odot i(a))=2\pi
q_2(i(a)\odot i(b)).\]

\begin{proposition}\label{prop.sumovertrees}
In the set-up of Theorem~\ref{thm.transfer}, if
$K\in\Hom^{-1}(V[1],V[1])$ satisfies the equation
$q_1K+Kq_1=\imath\pi-\Id_{V[1]}$, then the operators
$\imath^1_\infty$ and $\hat{Q}^1$ can be expressed as sums over
decorated rooted trees via the formulas
\[
\imath^1_\infty=\imath+\sum_{\Gamma\in
T_K}\frac{Z_\Gamma(\imath,\pi,K,q_i)}{|\operatorname{Aut}\Gamma|};\qquad
\hat{Q}^1_\infty=\delta_{[1]}+\sum_{\Gamma\in
T_\pi}\frac{Z_\Gamma(\imath,\pi,K,q_i)}{|\operatorname{Aut}\Gamma|}.
\]
In particular, for $n\geq 2$, the $n$-th higher bracket defining
the $L_\infty$ structure on $C[1]$ is
\[
\langle\,\rangle^{}_n=\sum_{\Gamma\in
T_{\pi,n}}\frac{Z_\Gamma(\imath,\pi,K,q_i)}{|\operatorname{Aut}\Gamma|},
\]
where ${\mathcal T}_{\pi,n}$ is the full subgroupoid of ${\mathcal
T}_{\pi}$ whose objects are rooted trees with exactly $n$ tail
vertices.
\end{proposition}
\begin{proof}
Let ${\mathcal V}_K$ be the groupoid whose objects are rooted
trees with a single internal vertex, of valence at least two. The
root edge is decorated by the operator $K$, the tail edges are
decorated by the operator $\Id_{V[1]}$ and the vertex is decorated
by the operator $q_r$, where $r$ is the number of tail edges.
denote by ${\mathcal V}_{K,n}$ the full subgroupoid  of ${\mathcal
V}_K$ whose object have exactly $n$ tail edges; note that the set
$V_{K,n}$ of isomorphism classes of objects of ${\mathcal
V}_{K,n}$ consists of a single element and that if $\Gamma$ is an
object in ${\mathcal V}_{K,n}$ with $n$-tail edges, then
\[
\frac{Z_\Gamma(\imath,\pi,K,q_i)}{|\operatorname{Aut}\Gamma|}=K
q_n.
\]
Also, let ${\mathcal I}$ be the groupoid whose objects are trees
consisting of a single directed edge decorated by the operator
$\imath$. Clearly, objects in ${\mathcal I}$ are all isomorphic to
each other and have no nontrivial automorphisms; operadic
evaluation gives, for an object $\Gamma$ in ${\mathcal I}$,
\[
\frac{Z_\Gamma(\imath,\pi,K,q_i)}{|\operatorname{Aut}\Gamma|}=\imath.
\]
  Each tree in
${\mathcal T}_K$ can be seen as the composition of exactly one
object of ${\mathcal V}_K$ and a forest consisting of trees in
${\mathcal T}_K$ and in ${\mathcal I}$. Hence, the usual
combinatorics of sums over graphs (see \cite{fiorenza:graphs} for
details) gives:
\begin{align*}
\imath+\sum_{\Gamma\in
T_K}\frac{Z_\Gamma(\imath,\pi,K,q_i)}
{|\operatorname{Aut}\Gamma|}&=\imath+\sum_{n\geq
2} \sum_{\Gamma_1\in
V_{K,n}}\!\!\frac{Z_{\Gamma_1}(\imath,\pi,K,q_i)}{|\operatorname{Aut}\Gamma_1|}\left(
\!\!\left(\imath+\sum_{\Gamma_2\in
T_K}\frac{Z_{\Gamma_2}(\imath,\pi,K,q_i)}
{|\operatorname{Aut}\Gamma_2|}\right)^{\!\odot
n}
\right)\\
&=\imath+\sum_{n\geq 2} K q_n\left(\left(\imath+
\sum_{\Gamma\in
T_K}\frac{Z_{\Gamma_2}(\imath,\pi,K,q_i)}
{|\operatorname{Aut}\Gamma|}\right)^{\!\odot
n}
\right).\\
\end{align*}
This shows that
\[
\imath+\sum_{\Gamma\in
T_K}\frac{Z_\Gamma(\imath,\pi,K,q_i)}{|\operatorname{Aut}\Gamma|}
\]
satisfies the same recursion as $\imath_\infty$; since this
recursion uniquely determines $\imath_\infty$, the formula for
$\imath_\infty$ is proved. The formula for $\hat{Q}^1$ is proved by
a completely similar argument.
\end{proof}

\bigskip
\section{The $L_\infty$ structure on $C_{\chi}$}
\label{sec.Linftyoncone}

By Quillen construction \cite{Qui}, the $L_{\infty}$ structure on
the differential graded Lie algebra $H_\chi$ is given by the
brackets
\[ {q}_k\colon \bigodot^k (H_\chi[1])\to H_\chi[1], \]
where ${q}_k=0$ for every $k\ge 3$,
\[ {q}_1(l,m(t,dt))=(-dl,-dm(t,dt))\]
and
\[q_2((l_1,m_1(t,dt))\odot(l_2,m_2(t,dt)
))=
(-1)^{\deg_{H_{\chi}}(l_1,m_1(t,dt))}([l_1,l_2],[m_1(t,dt),m_2(t,dt)]).\]

By Lemma~\ref{lem.definitionK} we can apply homotopy trasfer in
order to costruct an explicit $L_{\infty}$ structure on $C_{\chi}$
and an explicit $L_{\infty}$-morphism $\imath_{\infty}\colon
C_{\chi}\to H_{\chi}$ extending $\imath$.\\
According to Proposition~\ref{prop.sumovertrees}, the linear maps
of degree 1
\[ \langle\;\rangle_n^{}\colon \bigodot^n C_\chi[1]\to
C_\chi[1],\qquad n\ge 2,\]%
defining the induced $L_\infty$-algebra
structure on $C_\chi$ are explicitly described in terms of
summation over rooted trees. In our case, the  properties
\[
q_2(\image{K}\otimes\image{K})\subseteq\ker\pi\cap\ker K, \qquad
q_k=0\;\; \forall\; k\ge 3,\]%
imply that, fixing the number $n\ge 2$ of tails, there exists at
most one isomorphism class of trees  giving a nontrivial
contribution.
\begin{itemize}

\item{} $n=2$
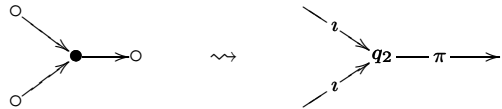
\[
\begin{xy}
,(-8,6)*{\circ};(0,0)*{\bullet}**\dir{-}?>*\dir{>}
,(-8,-6)*{\circ};(0,0)*{\bullet}**\dir{-}?>*\dir{>}
,(0,0)*{\bullet};(8,0)*{\circ}**\dir{-}?>*\dir{>}
\end{xy}
\qquad
\rightsquigarrow
\qquad
\begin{xy}
,(-10,6.66);(-6,4)*{\,\scriptstyle{\imath}\,}**\dir{-}
,(-10,-6.66);(-6,-4)*{\,\scriptstyle{\imath}\,}**\dir{-}
,(-6,4)*{\,\scriptstyle{\imath}\,};
(0,0)*{\,\,\scriptstyle{q_2}\,}**\dir{-}?>*\dir{>}
,(-6,-4)*{\,\scriptstyle{\imath}\,};(0,0)*{\,\,\scriptstyle{q_2}\,}**\dir{-}?>*\dir{>}
,(0,0)*{\,\,\scriptstyle{q_2}\,};
(8,0)*{\,\scriptstyle{\pi}\,}**\dir{-}
,(8,0)*{\,\scriptstyle{\pi}\,};
(16,0)**\dir{-}?>*\dir{>}
\end{xy}
\]
This graph gives
\[ \langle
{\gamma_1}\odot{\gamma_2}\rangle_2=\pi
q_2(\imath({\gamma_1})
\odot\imath({\gamma_2})).\]

\item{} $n>2$ :
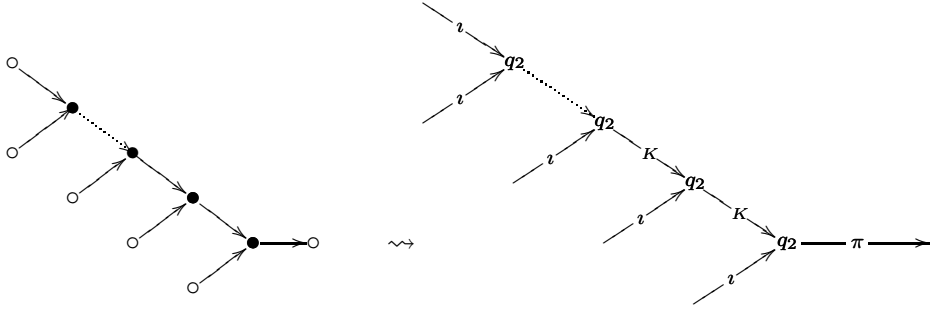
\[
\hskip -2.4em
\begin{xy}
,(-32,24)*{\circ};(-24,18)*{\bullet}**\dir{-}?>*\dir{>}
,(-32,12)*{\circ};(-24,18)**\dir{-}?>*\dir{>}
,(-24,6)*{\circ};(-16,12)*{\bullet}**\dir{-}?>*\dir{>}
,(-24,18)*{\circ};(-16,12)**\dir{.}?>*\dir{>}
,(-8,6)*{\circ};(0,0)*{\bullet}**\dir{-}?>*\dir{>}
,(-8,-6)*{\circ};(0,0)*{\bullet}**\dir{-}?>*\dir{>}
,(-16,12);(-8,6)*{\bullet}**\dir{-}?>*\dir{>}
,(-16,0)*{\circ};(-8,6)*{\bullet}**\dir{-}?>*\dir{>}
,(0,0)*{\bullet};(8,0)*{\circ}**\dir{-}?>*\dir{>}
\end{xy}
\qquad
\rightsquigarrow
\begin{xy}
,(-36,24)*{\,\,\scriptstyle{q_2}\,\,};(-24,16)*{\,\,\scriptstyle{q_2}\,\,}**\dir{.}?>*\dir{>}
,(-12,-8);(-7.2,-4.8)*{\,\scriptstyle{\imath}\,}**\dir{-}
,(-7.2,-4.8)*{\,\scriptstyle{\imath}\,};
(0,0)*{\,\,\scriptstyle{q_2}\,}**\dir{-}?>*\dir{>}
,(-24,0);(-19.2,3.2)*{\,\scriptstyle{\imath}\,}**\dir{-}
,(-19.2,3.2)*{\,\scriptstyle{\imath}\,};
(-12,8)*{\,\,\scriptstyle{q_2}\,\,}**\dir{-}?>*\dir{>}
,(-48,16);(-43.2,19.2)*{\,\scriptstyle{\imath}\,}**\dir{-}
,(-43.2,19.2)*{\,\scriptstyle{\imath}\,};
(-36,24)*{\,\,\scriptstyle{q_2}\,\,}**\dir{-}?>*\dir{>}
,(-48,32);(-43.2,28.8)*{\,\scriptstyle{\imath}\,}**\dir{-}
,(-43.2,28.8)*{\,\scriptstyle{\imath}\,};
(-36,24)*{\,\,\scriptstyle{q_2}\,\,}**\dir{-}?>*\dir{>}
,(-36,8);(-31.2,11.2)*{\,\scriptstyle{\imath}\,}**\dir{-}
,(-31.2,11.2)*{\,\scriptstyle{\imath}\,};
(-24,16)*{\,\,\scriptstyle{q_2}\,\,}**\dir{-}?>*\dir{>}
,(-12,8)*{\,\,\scriptstyle{q_2}\,\,};(-6,4)*{\,\scriptscriptstyle{K}\,}**\dir{-}
,(-6,4)*{\,\scriptscriptstyle{K}\,};
(0,0)*{\,\,\scriptstyle{q_2}\,}**\dir{-}?>*\dir{>}
,(-24,16)*{\,\,\scriptstyle{q_2}\,\,};(-18,12)*{\,\scriptscriptstyle{K}\,}**\dir{-}
,(-18,12)*{\,\scriptscriptstyle{K}\,};
(-12,8)*{\,\,\scriptstyle{q_2}\,\,}**\dir{-}?>*\dir{>}
,(0,0)*{\,\,\scriptstyle{q_2}\,};
(9.6,0)*{\,\scriptstyle{\pi}\,}**\dir{-}
,(9.6,0)*{\,\scriptstyle{\pi}\,};
(19.2,0)**\dir{-}?>*\dir{>}
\end{xy}
\]
\end{itemize}

This diagram gives, for every $n\ge 2$ the formula
\begin{multline*}
\langle{\gamma_1}\odot\cdots\odot
{\gamma_n}\rangle_n=\\
=\frac{1}{2}\sum_{\sigma\in S_n}\varepsilon(\sigma)\pi
{q}_2(\imath( {\gamma_{\sigma(1)}} \!)\odot
Kq_2(\imath({\gamma_{\sigma(2)}}) \odot\cdots\odot
Kq_2(\imath({\gamma_{\sigma(n-1)}})
\odot\imath({\gamma_{\sigma(n)}}))\cdots))\\
=\sum_{\substack{\sigma\in
S_n\\ \sigma(n-1)<\sigma(n)}}\varepsilon(\sigma)\pi
{q}_2(\imath( {\gamma_{\sigma(1)}} \!)\odot
Kq_2(\imath({\gamma_{\sigma(2)}})
\odot\cdots\odot Kq_2(\imath({\gamma_{\sigma(n-1)}})
\odot\imath({\gamma_{\sigma(n)}}))\cdots)).
\end{multline*}

\begin{remark}\label{rem.scalarextension}
The above construction of the $L_{\infty}$ structure on $C_{\chi}$
commutes
with tensor products of differential graded commutative
algebras.  This means that if $R$ is a DGCA, then the
$L_{\infty}$-algebra structure on the suspended mapping cone
of
$\chi\otimes{\rm id_R}\colon L\otimes R\to M\otimes R$ is naturally
isomorphic to the $L_{\infty}$-algebra
$C_{\chi}\otimes R$.
\end{remark}

A more refined description involving the original brackets in the
differential graded Lie algebras
$L$ and $M$ is obtained decomposing  the symmetric powers of
$C_\chi[1]$
into types:
\[
\bigodot^n \left(C_\chi[1]\right)=\bigodot^n
\Cone(\chi)=\bigoplus_{\lambda+\mu=n}
\left(\bigodot^\mu M\right)\otimes \left(\bigodot^\lambda
L[1]\right).
\]

The operation $\langle\,\rangle_2$ decomposes into
\begin{align*}
&{l_1}\odot {l_2}\mapsto
(-1)^{\deg_L^{}(l_1)}[l_1,l_2]\in L;\qquad
{m_1}\odot {m_2}\mapsto 0;\\
\\
&m\otimes l\mapsto
\dfrac{(-1)^{\deg_M^{}(m)+1}}{2}[m,\chi(l)]\in M
..
\end{align*}

For every $n\geq 2$ it is easy to see that
$\langle{\gamma_1}\odot\cdots\odot{\gamma_{n+1}}
\rangle_{n+1}$
can be nonzero only
if the multivector
${\gamma_1}\odot\cdots\odot{\gamma_{n+1}}$
belongs to
$\bigodot^{n}M\otimes L[1]$.
For $n\ge 2$, $m_1,\ldots,m_{n}\in M$ and $l\in L[1]$ the
formula for
$\langle\,\rangle_{n+1}$
described above becomes
\begin{multline*}
\langle {m_1}\odot\cdots\odot
{m_{n}}\otimes l\rangle_{n+1}=\\
=\sum_{\sigma\in S_{n}}\varepsilon(\sigma)\pi q_2((dt)
m_{\sigma(1)}\odot Kq_2((dt)m_{\sigma(2)}\odot\cdots\odot
Kq_2((dt)m_{\sigma(n)}\otimes t\chi(l))\cdots)).
\end{multline*}
Define recursively
a sequence of polynomials $\phi_i(t)\in \Q[t]\subseteq \mathbb{K}[t]$
and rational numbers
$I_n$ by the rule
\[ \phi_1(t)=t,\qquad I_n=\int_0^1\phi_n(t)dt,\qquad
\phi_{n+1}(t)=\int_0^t\phi_{n}(s)ds-tI_n.\]
By the definition of the homotopy operator $K$ we have, for every
$m\in M$
\[ K((\phi_n(t)dt)m)=\phi_{n+1}(t)m.\]
Therefore, for every  $m_1,m_2\in M$
   we have
\[
Kq_2((dt\cdot m_1)\odot
\phi_n(t){m_2})=-(-1)^{\deg_M^{}(m_1)}\phi_{n+1}(t) [m_1,m_2].\]
Therefore, we find:
\begin{multline*}
\langle m_1\odot\cdots\odot
m_{n}\otimes
l\,\rangle_{n+1}=\qquad\qquad\qquad\\
\quad=\sum_{\sigma\in S_{n}}\varepsilon(\sigma)\pi q_2((dt)
m_{\sigma(1)}\odot Kq_2((dt)m_{\sigma(2)} \odot\cdots\odot
Kq_2((dt)m_{\sigma(n)}\otimes
t\chi(l))\cdots))\\
   =(-1)^{1+\deg_M(m_{\sigma(n)})}\sum_{\sigma\in
S_{n}}\varepsilon(\sigma)\pi q_2((dt)
m_{\sigma(1)}\odot
Kq_2((dt)m_{\sigma(2)}
\odot\cdots\odot
\phi_2(t)[m_{\sigma(n)},\chi(l)]\cdots))\\
\qquad =(-1)^{n-1+\sum_{i=2}^n\deg_M(m_{\sigma(i)})}\sum_{\sigma\in
S_{n}}\varepsilon(\sigma)\pi q_2((dt)
m_{\sigma(1)}\odot\phi_{n}(t)[m_{\sigma(2)},\cdots,
[m_{\sigma(n)},\chi(l)]\cdots])\\
=(-1)^{n+\sum_{i=1}^n\deg_M(m_i)}I_n\sum_{\sigma\in
S_{n}}\varepsilon(\sigma)
[m_{\sigma(1)},[m_{\sigma(2)},\cdots,
[m_{\sigma(n)},\chi(l)]\cdots]]\in M.
\end{multline*}

\begin{theorem} For any morphism $\chi\colon L\to M$ of
differential graded Lie algebras, let
$\widetilde{C}(\chi)=(C_\chi, \hat{Q})$ be the
$L_\infty$-algebra structure defined on $C_\chi$ by the
above construction. Then
\[ \widetilde{C}\colon \mathbf{DGLA}^{\mathbf 2}\to
\mathbf{L}_\infty\]
is a functor
making the  diagram
\[\xymatrix{
\mathbf{DGLA}\ar[d]\ar[r]&\mathbf{L_\infty}\ar[d]\\
\mathbf{DGLA}^{\mathbf
2}\ar[r]\ar[ur]^{\widetilde{C}}&\mathbf{DG}}\]
commutative. If
$\sF\colon\mathbf{DGLA}^{\mathbf 2}\to \mathbf{L}_{\infty}$ has
the same properties, then for every morphism $\chi$ of
differential graded Lie algebras, the $L_{\infty}$-algebra
$\sF(\chi)$ is isomorphic to $\widetilde{C}(\chi)$.

\end{theorem}
\begin{proof}
The functoriality of $\widetilde{C}$ is clear:
in fact, for every commutative diagram
\[
\xymatrix{
      L_1 \ar[r]^{f_L}
\ar[d]_{\chi^{}_1} &
L_2 \ar[d]^{\chi^{}_2}\\
    M_1  \ar[r]^{f_M} &
M_2\\
    }
\]
of  morphisms of differential graded Lie algebras, the natural map
$(f_L,f_M)\colon \widetilde{C}({\chi^{}_1})\to \widetilde{C}({\chi^{}_2})$ is a linear
$L_{\infty}$-morphism.\\
If $\sF\colon\mathbf{DGLA}^{\mathbf 2}\to \mathbf{L}_{\infty}$ has
the same properties as $\widetilde{C}$, according to 
Proposition~\ref{prop.esistequasiiso}
there exists two injective $L_{\infty}$ quasi-isomorphisms 
\[ \imath_{\infty}\colon \widetilde{C}(\chi)\to H_{\chi},\qquad
\hat{\imath}_{\infty}\colon \sF(\chi)\to H_{\chi}\]
with the same linear term $\imath$. The composition of $\imath_{\infty}$ with a 
left inverse of $\hat{\imath}_{\infty}$ is an isomorphism of $L_{\infty}$-algebras.
\end{proof}

\begin{remark}
As an instance of functoriality, note that the projection on
the first factor $p_1\colon \widetilde{C}(\chi)\to L$ is a
linear morphism of $L_\infty$-algebras. To see this, 
consider the morphism in $\mathbf{DGLA}^{\mathbf 2}$
\[
\xymatrix{
      L \ar[r]^{\Id_L}
\ar[d]_{\chi} &
L \ar[d]^{}\\
    M  \ar[r]^{} &
0\\
    }
\]
\end{remark}

We also have an explicit expression for the coefficients
$I_n$ appearing in the formula for $\langle\,\rangle_{n+1}$;
in the next lemma we show that they are, up to a sign, the
Bernoulli numbers. 
\begin{lemma}
For every $n\geq1$ we have $I_n=-B_n/n!$, 
where the $B_n$ are the Bernoulli numbers, i.e., the rational
numbers defined by the series expansion identity
\[
\sum_{n=0}^{\infty}B_n\frac{x^n}{n!}=\frac{x}{e^x-1}
=1-\frac{x}{2}+\frac{x^2}{12}-\frac{x^4}{720}+\frac{x^6}{30240}-
\frac{x^8}{1209600}+\cdots
\]
\end{lemma}

\begin{proof} Keeping in mind the definition of $B_n$,
we have to prove that
\[
1-\sum_{n=1}I_n
x^n=
\frac{x}{e^{x}-1}.
\]
Consider the polynomials $\psi_0(t)=1$ and
$\psi_n(t)=\phi_n(t)-I_n$ for $n\geq 1$. Then, for any $n\geq 1$,
\[\qquad\displaystyle\frac{d}{dt}\psi_n(t)=\psi_{n-1}(t),\qquad
\int_0^1\psi_n(t)dt=0.\]%
Setting
\[
F(t,x)=\sum_{n=0}^\infty \psi_n(t)x^n,
\]
we have
\[
\frac{d}{dt}F(t,x)=\sum_{n=1}^\infty\psi_{n-1}(t)x^n=xF(t,x),\qquad
\int_0^1F(t,x)dt=1.
\]
Therefore,
\[
F(t,x)=F(0,x)e^{tx},
\]
\[
1=\int_0^1F(t,x)dt=F(0,x)\int_0^1e^{tx}dt=
F(0,x)\frac{e^{x}-1}{x},
\]
and then
\[
F(0,x)=\frac{x}{e^{x}-1}.
\]
Since $\psi_n(0)=-I_n$ for any $n\geq 1$ we get
\[
\frac{x}{e^{x}-1}=F(0,x)=1-\sum_{n=1}^\infty
I_nx^n.
\]
In fact an alternative proof of the equality $I_n=-B_n/n!$ can be
done by observing that  the polynomials $n!\psi_n(t)$ satisfy the
recursive relations of the Bernoulli polynomials (see e.g.
\cite{remmert}).
\end{proof}

Summing up the results of this Section,
we have the following explicit description of the
$L_\infty$ algebra $\widetilde{C}(\chi)$.

\begin{theorem}\label{thm.coefficienti}
The $L_\infty$ algebra $\widetilde{C}(\chi)$ is defined by
the multilinear maps
\[
\langle\,\rangle_n\colon \bigodot^nC_\chi[1]\to C_\chi[1],
\]
given by
\[
\langle (l,m)\rangle_1^{}=(-dl, -\chi(l)+dm);
\]
\[
\langle l_1\odot l_2\rangle_2^{}=(-1)^{\deg_L(l_1)}[l_1,l_2];
\]
\[
\langle m\otimes
l\rangle_2^{}=\frac{(-1)^{\deg_M(m)+1}}{2}[m,\chi(l)];
\]
\[
\langle m_1\odot
m_2\rangle_2^{}=0;
\]
\[
\langle{m_1}\odot\cdots\odot
{m_{n}}\otimes
l_1\odot\cdots\odot l_k\rangle_{n+k}^{}= 0, \qquad n+k\geq 3,
k\neq 1;
\]
and
\begin{multline*}
\langle {m_1}\odot\cdots\odot
{m_{n}}\otimes
l\,\rangle_{n+1}=\\
=-(-1)^{\sum_{i=1}^n\deg_M(m_i)}\frac{B_n}{n!}\sum_{\sigma\in
S_{n}}\varepsilon(\sigma)
[m_{\sigma(1)},[m_{\sigma(2)},\cdots,
[m_{\sigma(n)},\chi(l)]\cdots]],\qquad n\geq 2;
\end{multline*}
where the $B_n$ are the Bernoulli numbers, i.e., the rational
numbers defined by the series expansion identity
\[
\sum_{n=0}^{\infty}B_n\frac{x^n}{n!}=\frac{x}{e^x-1}
=1-\frac{x}{2}+\frac{x^2}{12}-\frac{x^4}{720}+\frac{x^6}{30240}-
\frac{x^8}{1209600}+\cdots
\]

\end{theorem}

\begin{remark}
Via the decalage isomorphism
$\odot^n(C_\chi[1])\xrightarrow{\sim}(\wedge^n C_\chi)[n]$, the
linear maps $\langle\,\rangle_n$ defining the $L_\infty$-algebra
$\widetilde{C}(\chi)$ correspond to multilinear
operations $[\,]_n^{}\colon\wedge^n C_\chi\to C_\chi [2-n]$ on
$C_\chi$. In particular, the linear map
$\langle\,\rangle^{}_1$ corresponds to the differential
$\delta$ on $C_\chi$
\[
\delta\colon(l,m)\mapsto(dl,\chi(l)-dm),
\]
whereas the map $\langle\,\rangle_2^{}$ corresponds to the
following degree zero bracket
\[
[\,]_2^{}\colon C_\chi\wedge C_\chi \to C_\chi
\]  
\[
[l_1,l_2]_2^{} = [l_1,l_2];\qquad [m,l]_2^{}=
\frac{1}{2}[m,\chi(l)];\qquad [m_1,m_2]_2^{}=
0.
\]
Note that this is precisely the naive bracket described in
the introduction.
\end{remark}

\begin{remark} The occurrence of Bernoulli numbers is not
surprising: it had already been noticed by K.~T.~Chen
\cite{chen} how Bernoulli numbers are related to the
coefficients of the Baker-Campbell-Hausdorff formula.\\
More
recently, the relevance of Bernoulli numbers in
deformation theory  has been also
remarked by Ziv Ran in versions v3-v6 of \cite{ranATOMS}. 
In particular, Ziv Ran's
``JacoBer" complex seems to
be closely related with the coderivation $\hat{Q}$ defining the
$L_\infty$ structure on $C_\chi$.\\
Bernoulli numbers also appear in some expressions of the gauge
equivalence in a differential graded Lie algebra
\cite{sullivan,getzler}. In fact the relation $x=e^a\ast y$
   can be written as
\[ x-y=\frac{e^{\ad_{a}}-1}{\ad_{a}}([a,y]-da).\]
Applying to both sides the inverse of the
operator $\dfrac{e^{\ad_{a}}-1}{\ad_{a}}$
we get
\[ da=[a,y]-\sum_{n\ge 0}\frac{B_n}{n!}\ad_{a}^n(x-y).\]
\\
The multilinear brackets $\langle\;\rangle_n$
on $\Cone(\chi)=C_\chi[1]$ can be related to the Koszul (or `higher
derived') brackets $\Phi_n$ of a differential graded Lie
algebra as follows. Let $(M,\de, [\,,\,])$ be a differential
graded Lie algebra; the Koszul brackets
\[ \Phi_n\colon \bigodot^n M\to M,\qquad n\ge 1\]
are the degree 1 linear maps defined as $\Phi_1=0$ and for $n\ge 2$
\[ \Phi_n(m_1\cdots m_n)=\frac{1}{n!}\sum_{\sigma\in
\Sigma_n} \varepsilon(\sigma) [\cdot[[\de
m_{\sigma(1)},m_{\sigma(2)}],m_{\sigma(3)}],\cdots,
m_{\sigma(n)}].\]
Let $L$ be the differential graded Lie subalgebra of $M$
given by $L:=\de M$ and let $\chi\colon L\to M$ be the
inclusion. We can identify $M$ with the image of the
injective linear map
$M\hookrightarrow
\Cone(\chi)$ given by $m\mapsto (\de m, m)$.
Then we have
$\langle (\de m,m)\rangle_1=0$,
\[
\langle (\de m_1,m_1)\odot (\de
m_2,m_2)\rangle^{}_2=(\de\Phi_2(m_1,m_2),
\Phi_2(m_1,m_2))
\]
and, for $n\geq 2$,
\[
\langle
(\de m_{1},m_{1})\odot\cdots \odot(\de
m_{{n+1}},m_{{n+1}})
\rangle_{n+1}^{}=(0, B_n(-1)^{n}(n+1)
\Phi_{n+1}(m_{1}\odot\cdots \odot m_{{n+1}})).
\]
Since the
multilinear operations $\langle\,\rangle_n$ define an
$L_\infty$-algebra structure on $C_\chi=\Cone(\chi)[-1]$,
    they satisfy a sequence of quadratic
relations. Due to the
above mentioned correspondence with the Koszul brackets,
these relations are translated into a sequence of
differential/quadratic relations between the
odd Koszul brackets, defined as
$\{m\}^{}_1=0$ and
\[ \{m_1,\cdots, m_n\}^{}_n=\frac{1}{n!}\sum_{\sigma\in
\Sigma_n} \varepsilon(\sigma) (-1)^\sigma[\cdot[[\de
m_{\sigma(1)},m_{\sigma(2)}],m_{\sigma(3)}],\cdots,
m_{\sigma(n)}]\]
for $n\geq 2$. For
instance, if
$m_1,m_2,m_3$ are homogeneous elements of degree
$i_1,i_2,i_3$ respectively, then
\begin{align*}
\{\{m_1,m_2\}_2^{},m_3\}^{}+&(-1)^{i_1i_2+i_1i_3}
\{\{m_2,m_3\}_2^{},m_1\}_2^{}\\
&+(-1)^{i_2i_3+i_1i_3}\{\{m_3,m_1\}_2^{},m_2\}^{}_2
=3/2\cdot\de \{m_1,m_2,m_3\}_3^{}.
\end{align*}
The occurrence of Bernoulli numbers in the $L_\infty$-type structure
defined by the higher Koszul brackets has been recently remarked by
K.~Bering in \cite{bering}.
\end{remark}

\bigskip
\section{The Maurer-Cartan functor}
\label{sec.maurercartan}

Having introduced an $L_\infty$ structure on $C_\chi$ in
Section~\ref{sec.Linftyoncone}, we have a corresponding
Maurer-Cartan functor \cite{fuka,K}
$\MC_{C_\chi}\colon\mathbf{Art}\to\mathbf{Set}$, defined as
\[\MC_{C_\chi}(A)=\left\{ \gamma\in
C_\chi[1]^0\otimes\mathfrak{m}_A
\;\strut\left\vert\;
\sum_{n\ge1}\frac{\langle\gamma^{\odot
n}\rangle_n}{n!}=0\right.\right\},\qquad
A\in \mathbf{Art}.
\]

Writing $\gamma=(l,m)$, with $l\in
L^1\otimes\mathfrak{m}_A$ and $m\in
M^0\otimes\mathfrak{m}_A$, the Maurer-Cartan equation becomes
\begin{align*}
0&=\sum_{n=1}^\infty\frac{\langle(l,m)^{\odot
n}\rangle_n}{n!}\\
&=
\langle (l,m)\rangle_1+\frac{1}{2}\langle
l^{\odot 2}\rangle_2+\langle
m\otimes l\rangle_2+\frac{1}{2}\langle
m^{\odot 2}\rangle_2+\sum_{n\geq 2}
\frac{n+1}{(n+1)!}\langle m^{\odot n}\otimes
l\rangle_{n+1}\\
&=
\left(-dl-\frac{1}{2}[l,l],
-\chi(l)+dm -\frac{1}{2}[m,\chi(l)]
+\sum_{n\geq 2}
\frac{1}{n!}\langle m^{\odot n}\otimes
l\rangle_{n+1}\right)\in(L^2\oplus M^1)\otimes\mathfrak{m}_A.
\end{align*}

According to Theorem~\ref{thm.coefficienti}, since
$\deg_M(m)=\deg_{C_{\chi[1]}}(m)=0$, we have
\[ \langle m^{\odot n}\otimes
l\rangle_{n+1}=-\frac{B_n}{n!} \sum_{\sigma\in S_{n}}[m,[m,\cdots,
[m,\chi(l)]\cdots]]=-B_n\ad_m^n(\chi(l)),\]%
where for $a\in M^0\otimes \mathfrak{m}_A$ we denote by
$\ad_a\colon M\otimes\mathfrak{m}_A\to M\otimes\mathfrak{m}_A$ the
operator $\ad_a(y)=[a,y]$.

The Maurer-Cartan equation on $C_{\chi}$ is therefore  equivalent
to
\[
\begin{cases}
\displaystyle{dl+\frac{1}{2}[l,l]=0}\\
\\
\displaystyle{\chi(l)-dm+\frac{1}{2}[m,\chi(l)]+
\sum_{n=2}^{\infty}\frac{B_n}{n!} \ad_m^n(\chi(l))=0}.
\end{cases}
\]
Since $B_0=1$ and $B_1=-\dfrac{1}{2}$, we can write the
second equation as
\begin{multline*}
0=\chi(l)-dm+\frac{1}{2}[m,\chi(l)]+
\sum_{n=2}^{\infty}\frac{B_n}{n!} \ad_m^n(\chi(l))\\
=[m,\chi(l)]-dm+
\sum_{n=0}^{\infty}\frac{B_n}{n!}\ad_m^n(\chi(l))
=[m,\chi(l)]-dm+\frac{\ad_m}{e^{\ad_m}-1}(\chi(l)).
\end{multline*}%
Applying the invertible operator $\dfrac{e^{\ad_m}-1}{\ad_m}$  we get
\[
0=\chi(l)+\dfrac{e^{\ad_m}-1}{\ad_m}([m,\chi(l)]-dm).\]

In the right side of the last formula we recognize the explicit
description of the gauge action
\[ \exp(M^0\otimes \mathfrak{m}_A)\times M^1\otimes
\mathfrak{m}_A\mapor{\ast} M^1\otimes \mathfrak{m}_A,\]
\[ e^a\ast
y=y+\sum_{n=0}^{+\infty}\frac{\ad^n_a}{(n+1)!}([a,y]-da)
=y+\frac{e^{\ad_a}-1}{\ad_a}([a,y]-da).\]%

Therefore, the
Maurer-Cartan equation for the $L_\infty$-algebra
structure on $C_\chi$ is equivalent to
\[
\begin{cases}
dl+\dfrac{1}{2}[l,l]=0\\
e^m\ast\chi(l)=0.
\end{cases}
\]

\bigskip
\section{Homotopy equivalence and the deformation functor}
\label{sec.gaugeequivalence}

Recall that the deformation functor associated to an
$L_\infty$-algebra $\mathfrak{g}$ is
$\Def_{\mathfrak{g}}=\MC_{\mathfrak{g}}/\sim$, where $\sim$
denotes homotopy equivalence of solutions of the Maurer-Cartan
equation: two elements $\gamma_0$ and $\gamma_1$ of
$\MC_{\mathfrak g}(A)$ are called homotopy equivalent if there
exists an element $\gamma(t,dt)\in\MC_{{\mathfrak g}[t,dt]}(A)$
with $\gamma(0)=\gamma_0$ and $\gamma(1)=\gamma_1$. It is possible
to prove that homotopy equivalence is an equivalence relation; in
this paper we do not need this fact.

We have already  described  the functor $\MC_{C_\chi}$ in terms of
the Maurer-Cartan equation in $L$ and the gauge action in $M$. Now
we want to prove a similar result for   the homotopy equivalence
on $\MC_{C_\chi}$. We need some preliminary results.

\begin{proposition}\label{prop.hvg}
Let $(L,d,[~,~])$ be a differential graded Lie algebra such that:
\begin{enumerate}

\item $L=M\oplus C\oplus D$ as graded vector spaces.

\item $M$ is a differential graded subalgebra of $L$.

\item  $d\colon C\to D[1]$ is an isomorphism of graded vector
spaces.
\end{enumerate}
Then, for every $A\in \mathbf{Art}$ there exists a bijection
\[ \alpha\colon \MC_M(A)\times (C^0\otimes
\mathfrak{m}_A)\mapor{\sim}\MC_L(A),\qquad (x,c)\mapsto e^c\ast
x.\]
\end{proposition}

\begin{proof} This is essentially proved in \cite[Section 5]{SchSta}
using induction on the length of $A$ and the
Baker-Campbell-Hausdorff formula.  Here we sketch a different
proof based on formal
theory of deformation functors \cite{Sch,Rim,FM1,ManettiDGLA}.\\
The map $\alpha$ is a natural transformation of homogeneous
functors, so it is sufficient to show that $\alpha$ is bijective
on tangent spaces and injective on obstruction spaces. Recall that
the tangent space of $\MC_L$ is $Z^1(L)$, while its obstruction
space is  $H^2(L)$. The functor $A\mapsto
C^0\otimes\mathfrak{m}_A$ is smooth with tangent space $C^0$ and
therefore tangent and obstruction spaces of the functor
\[ A\mapsto \MC_M(A)\times(C^0\otimes\mathfrak{m}_A)\]
are respectively $Z^1(M)\oplus C^0$ and $H^2(M)$. The tangent map
is
\[ Z^1(M)\oplus C^0\ni (x,c)\mapsto e^c\ast x=x-dc\in
Z^1(M)\oplus d(C^0)=Z^1(M)\oplus D^1=Z^1(L)\]%
and it is an isomorphism. The inclusion $M\hookrightarrow L$ is a
quasiisomorphism, therefore the obstruction to lifting $x$ in $M$
is equal to the obstruction to lifting $x=e^0\ast x$ in $L$. We
conclude the proof by observing that, according to \cite[Prop.
7.5]{FM1}, \cite[Lemma 2.20]{ManettiDGLA}, the obstruction maps of
Maurer-Cartan functor are invariant under the gauge action.
\end{proof}

\begin{corollary}\label{cor.homvsgauge1}
Let $M$ be a differential graded Lie algebra, $L=M[t,dt]$ and
$C\subseteq M[t]$ the subspace consisting of polynomials $g(t)$
with $g(0)=0$. Then for every $A\in \mathbf{Art}$ the map
$(x,g[t])\mapsto e^{g(t)}\ast x$ induces an isomorphism
\[\MC_M(A)\times (C^0\otimes \mathfrak{m}_A)\simeq\MC_L(A).\]
\end{corollary}
\begin{proof} The data $M,C$ and $D=d(C)$ satisfy the condition of
Proposition~\ref{prop.hvg}.\end{proof}

\begin{corollary}
Let $M$ be a differential graded Lie algebra.
Two elements $x_0,x_1\in \MC_M(A)$ are gauge equivalent if and
only if they are homotopy equivalent.
\end{corollary}

\begin{proof}
If $x_0$ and $x_1$ are gauge equivalent, then there exists $g\in
M^0\otimes\mathfrak{m}_A$ such that $e^g\ast x_0=x_1$. Then, by
Corollary~\ref{cor.homvsgauge1}. $x(t)=e^{t\, g}\ast x_0$ is an
element of $\MC_{M[t,dt]}(A)$ with $x(0)=x_0$ and $x(1)=x_1$,
i.e.,
$x_0$ and $x_1$ are homotopy equivalent.\\
Vice versa, if $x_0$ and $x_1$ are homotopy equivalent, there
exists $x(t)\in \MC_{M[t,dt]}(A)$ such that $x(0)=x_0$ and
$x(1)=x_1$. By Corollary~\ref{cor.homvsgauge1}., there exists
$g(t)\in M^0[t]\otimes \mathfrak{m}_A$ with $g(0)=0$ such that
$x(t)=e^{g(t)}\ast x_0$. Then $x_1=e^{g(1)}\ast x_0$, i.e., $x_0$
and $x_1$ are gauge equivalent.
\end{proof}

\begin{theorem}
Let $\chi\colon L\to M$ be a morphism of differential graded Lie
algebras and let $(l_0,m_0)$ and $(l_1,m_1)$ be elements of
$\MC_{C_\chi}(A)$. Then $(l_0,m_0)$ is homotopically equivalent to
$(l_1,m_1)$ if and only if there exists $(a,b)\in
C_{\chi}^0\otimes \mathfrak{m}_A$ such that
\[ l_1=e^a\ast l_0,\qquad e^{m_1}=e^{db}e^{m_0}e^{-\chi(a)}.\]
\end{theorem}

\begin{remark} The condition $e^{m_1}=e^{db}e^{m_0}e^{-\chi(a)}$
can be also written as $m_1\bullet \chi(a)=db\bullet m_0$, where
$\bullet$ is the Baker-Campbell-Hausdorff product in the nilpotent
Lie algebra $M^0\otimes\mathfrak{m}_A$.\\
As a consequence, we get that in this case the homotopy
equivalence is induced by a group action; this is  false for
general $L_{\infty}$-algebras.
\end{remark}

\begin{proof} We shall say that two elements
$(l_0,m_0)$, $(l_1,m_1)$  are gauge equivalent if and only if
there exists $(a,b)\in C_{\chi}^0\otimes \mathfrak{m}_A$ such that
\[ l_1=e^a\ast l_0,\qquad e^{m_1}=e^{db}e^{m_0}e^{-\chi(a)}.\]%
We first show that homotopy implies gauge. Let $(l_0,m_0)$ and
$(l_1,m_1)$ be homotopy equivalent elements of $\MC_{C_\chi}(A)$.
Then there exists an element $(\tilde{l},\tilde{m})$  of
$\MC_{C_\chi[s,ds]}(A)$ with
$(\tilde{l}(0),\tilde{m}(0))=(l_0,m_0)$ and
$(\tilde{l}(1),\tilde{m}(1))=(l_1,m_1)$. According to
Remark~\ref{rem.scalarextension}, the Maurer-Cartan equation for
$(\tilde{l},\tilde{m})$ is
\[
\begin{cases}
d\tilde{l}+\dfrac{1}{2}[\tilde{l},\tilde{l}]=0\\
e^{\tilde{m}}\ast\chi(\tilde{l})=0
\end{cases}
\]
The first of the two equations above tells us that $\tilde{l}$ is
a solution of the Maurer-Cartan equation for $L[s,ds]$. So, by
Corollary~\ref{cor.homvsgauge1}, there exists a degree zero
element $\lambda(s)$ in $L[s]\otimes\mathfrak{m}_A$ with
$\lambda(0)=0$ such that $\tilde{l}=e^{\lambda}\ast l_0$.
Evaluating at $s=1$ we find $l_1=e^{\lambda_1}\ast l_0$. As a
consequence of $\tilde{l}=e^{\lambda}\ast l_0$, we also have
$\chi(\tilde{l})=e^{\chi(\lambda)}\ast\chi(l_0)$. Set
$\tilde{\mu}=\tilde{m}\bullet\chi(\lambda)\bullet m_0$, so that
$\tilde{m}=\tilde{\mu}\bullet m_0\bullet(-\chi(\lambda))$ and the
second Maurer-Cartan equation is reduced to
$e^{\tilde{\mu}}\ast(e^{m_0}\ast\chi(l_0))=0$, i.e., to
$e^{\tilde{\mu}}\ast 0=0$, where we have used the fact that
$(l_0,m_0)$ is a solution of the Maurer-Cartan equation in
$C_\chi$. This last equation is equivalent to the equation
$d\tilde{\mu}=0$ in $C_\chi[s,ds]\otimes\mathfrak{m}_A$. If we
write $\tilde{\mu}(s,ds)=\mu^0(s)+ds\,\mu^{-1}(s)$, then the
equation $d\tilde{\mu}=0$ becomes
\[
\begin{cases}
\dot{\mu}^0-d^{}_M\mu^{-1}=0\\
d^{}_M\mu^0=0,
\end{cases}
\]
where $d^{}_M$ is the differential in the DGLA $M$.
The solution is, for any fixed $\mu^{-1}$,
\[
\mu^0(s)=\int_0^s d\sigma\, d^{}_M\mu^{-1}(\sigma)=-d^{}_M
\int_0^s d\sigma\,
\mu^{-1}(\sigma)
\]
Set $\nu=-\int_0^1 d s\,
\mu^{-1}(s)$. Then $m_1=\tilde{m}(1)=(d_M\nu)\bullet
m_0\bullet (-\chi(\lambda_1))$. Summing up, if $(l_0,m_0)$
and
$(m_1,l_1)$ are homotopy equivalent, then there exists
$(d\nu,\lambda_1)\in (d M^{-1}\otimes\mathfrak{m}_A)\times
(L^0\otimes\mathfrak{m}_A)$ such that
\[
\begin{cases}
l_1=e^{\lambda_1}\ast l_0\\
m_1=d\nu\bullet m_0\bullet (-\chi(\lambda_1)),
\end{cases}
\]
i.e., $(l_0,m_0)$
and
$(m_1,l_1)$ are gauge equivalent.

We now show that gauge
implies homotopy. Assume  $(l_0,m_0)$
and
$(m_1,l_1)$ are gauge equivalent. Then there exists
$(d\nu,\lambda_1)\in (d M^{-1}\otimes{\mathfrak m})\times
(L^0\otimes{\mathfrak m})$ such that
\[
\begin{cases}
l_1=e^{\lambda_1}*l_0\\
m_1=d\nu\bullet m_0\bullet (-\chi(\lambda_1)).
\end{cases}
\]
Set $\tilde{l}(s,ds)=e^{s\lambda_1}*l_0$. By
Corollary~\ref{cor.homvsgauge1}, $\tilde{l}$ satisfies the
equation $d\tilde{l}+\frac{1}{2}[\tilde{l},\tilde{l}]=0$. Set
$\tilde{m}=(d(s\nu))\bullet m_0\bullet(-\chi(s\lambda_1))$.
Reasoning as above, we find
\[
e^{\tilde m}*\chi(\tilde{l})=e^{d(s\nu)}*0=0.
\]
Therefore, $(\tilde{l},\tilde{m})$ is a solution of the
Maurer-Cartan equation in $C_\chi[s,ds]$. Moreover
$\tilde{l}(0)=l_0$, $\tilde{l}(1)=l_1$, $\tilde{m}(0)=m_0$ and
$\tilde{m}(1)=d\nu\bullet m_0\bullet (-\chi(\lambda_1))=m_1$, i.e.
$(l_0,m_0)$ and $(m_1,l_1)$ are homotopy equivalent.\end{proof}

\bigskip
\section{Examples and applications}
\label{sec.applicationtodeftheory}

Let $\chi\colon L\to M$ be a morphism of differential graded Lie
algebras over a field $\K$ of characteristic 0. In the paper
\cite{semireg} one of the authors has introduced, having in mind
the example of embedded deformations, the notion of Maurer-Cartan
equation and gauge action for the triple $(L,M,\chi)$; these
notions reduce to the standard Maurer-Cartan equation and gauge
action of $L$ when $M=0$. More precisely there are defined two
functors of Artin rings $\MC_{\chi},\Def_{\chi}\colon
\mathbf{Art}\to \mathbf{Set}$, in the following way:
\[ \MC_{\chi}(A)=
\left\{(x,e^a)\in (L^1\otimes\mathfrak{m}_A)\times \exp(M^0\otimes
\mathfrak{m}_A)\mid
dx+\frac{1}{2}[x,x]=0,\;e^a\ast\chi(x)=0\right\},\]
\[ \Def_\chi(A)=\frac{\MC_{\chi}(A)}{\text{gauge equivalence}},\]
where two solutions of the Maurer-Cartan equation are gauge
equivalent if they belong to the same orbit of the gauge action
\[(\exp(L^0\otimes\mathfrak{m}_A)\times
\exp(dM^{-1}\otimes\mathfrak{m}_A))
\times \MC_{\chi}(A)\mapor{\ast}\MC_{\chi}(A)\]%
given by the formula
\[ (e^l, e^{dm})\ast (x,e^a)=(e^l\ast x,
e^{dm}e^ae^{-\chi(l)})=(e^l\ast x,
e^{dm\bullet a\bullet (-\chi(l))}).\]%
The computations of Sections \ref{sec.maurercartan} and
\ref{sec.gaugeequivalence} show that $\MC_{\chi}$ and
$\Def_{\chi}$ are canonically isomorphic to the functors
$\MC_{\widetilde{C}(\chi)}$ and $\Def_{\widetilde{C}(\chi)}$
associated with the $L_{\infty}$ structure on $C_{\chi}$.

\begin{example}\label{exa.hilbertfunctor}
Let $X$ be a compact complex manifold and let $Z\subset X$ be a
smooth subvariety. Denote by $\Theta_X$ the  holomorphic tangent
sheaf of $X$ and by $N_{Z|X}$ the normal sheaf of $Z$ in $X$.\\
Consider the short exact sequence of complexes
\[0\to \ker{\pi}\mapor{\chi}A^{0,*}_X(\Theta_X)
\mapor{\pi}A^{0,*}_Z(N_{Z|X})\to 0.\]%
It is proved in \cite{semireg} that there exists a natural
isomorphism between  $\Def_{\chi}$ and the functor of embedded
deformations of $Z$ in $X$. Therefore, the $L_\infty$ algebra
$\widetilde{C}(\chi)$ governs embedded deformations in this
case.\\
Note  that the DGLA $A^{0,*}_Z(\Theta_Z)$ governs the deformations
of $Z$ and that the natural transformation
\[ \Def_{\widetilde{C}(\chi)}=\Def_{\chi}\to \Def_{A^{0,*}_Z(\Theta_Z)},\]
\[\{\text{Embedded deformations of }Z\}\to
\{\text{Deformations of }Z\}\]%
is induced by the morphism in $\mathbf{DGLA^2}$ given  by the
diagram
\[ \begin{matrix}
\ker{\pi}&\to&A^{0,*}_Z(\Theta_Z)\\
\mapver{\chi}&&\mapver{}\\
A^{0,*}_X(\Theta_X)&\to&0.
\end{matrix}\]

\end{example}

The next result  was proved in \cite{semireg} using the theory of
extended deformation functors; here we can prove it in a more
standard way.

\begin{theorem}\label{thm.basic}
Consider a commutative diagram of morphisms of differential graded
Lie algebras
\[
\xymatrix{
      L_1 \ar[r]^{f_L}
\ar[d]_{\chi^{}_1} &
L_2 \ar[d]^{\chi^{}_2}\\
    M_1  \ar[r]^{f_M} &
M_2\\
    }
\]
and assume that $(f_L,f_M)\colon C_{\chi^{}_1}\to C_{\chi^{}_2}$
is a quasiisomorphism of complexes (e.g. if both $f_L$ and $f_M$
are quasiisomorphisms). Then the natural transformation
$\Def_{\chi^{}_1}\to \Def_{\chi^{}_2}$ is an isomorphism.
\end{theorem}

\begin{proof} The map $(f_L,f_M)\colon \widetilde{C}(\chi^{}_1)\to \widetilde{C}(\chi^{}_2)$
is a linear quasi-isomorphism of $L_{\infty}$-algebras and then 
induces an isomorphism of the associated
deformation functors \cite{K}.
\end{proof}

\begin{example}
Let $\pi\colon A\to B$ be a surjective morphism of associative
$\K$-algebras and denote by $I$ its kernel. The algebra $B$ is an
$A$-module via $\pi$; this makes $B$ a trivial $I$-module. Let $K$
be the suspended Hochschild complex
\[
K=\Hoch^\bullet(I,B)[-1],
\]
Note that the differential $d$ of $K$ is identically zero if and
only if $I\!\cdot\! I=0$.

The natural map
\[ \alpha\colon\Hoch^\bullet(A,A)\to
K[1]=\Hoch^\bullet(I,B)\] is a surjective morphism of complexes,
and its kernel
\[ \ker\alpha=\{f\mid f(I^{\otimes})\subseteq I\}\]
is a Lie subalgebra of $\Hoch^\bullet(A,A)$ endowed with the
Hochschild bracket. Denote by $\chi\colon
\ker\alpha\hookrightarrow \Hoch^\bullet(A,A)$ the inclusion. Since
$\chi$ is injective, the projection on the second factor induces a
quasiisomorphism of differential complexes
\[
{\rm pr}_2\colon C_\chi\to \coker(\chi)[-1]\simeq K,
\]
where the isomorphism on the right is induced by the map $\alpha$.
Therefore we have a canonical $L_{\infty}$ structure (defined up
to homotopy) on $K$. This gives a Lie structure on the cohomology
of $K$, which is not trivial in general: consider for instance the
exact sequence

\[ 0\to \K\varepsilon\to \frac{\K[\varepsilon]}{(\varepsilon^2)}\mapor{\pi}
\K\to 0\] and
\[f\in K^1=H^1(K),\qquad f(\varepsilon)=1.\]%
Choose as a lifting the linear map $g\colon
\K[\varepsilon]/(\varepsilon^2)\to
\K[\varepsilon]/(\varepsilon^2)$,
\[ g(1)=0,\quad g(\varepsilon)=1.\]%
Then
\[ dg(\varepsilon\otimes \varepsilon)=2\varepsilon\]
and so $dg\in\ker\alpha$. Therefore, $(dg,g)$ is a closed element
of $C_\chi^1$ representing the cohomology class $f\in H^1(K)$ and
so
\[
[f,f]=\alpha({\rm pr}_2([(dg,g),(dg,g)]_2{}))=\alpha([g,dg]).
\]
One computes
\begin{align*}
[f,f](\varepsilon\otimes\varepsilon)&=\pi([g,dg](\varepsilon\otimes\varepsilon))\\
&=\pi(g(dg
(\varepsilon\otimes\varepsilon))
-dg(g(\varepsilon)\otimes\varepsilon)+dg(\varepsilon\otimes
g(\varepsilon)))\\
&=\pi(g(2\varepsilon)-dg(1\otimes\varepsilon)+dg(\varepsilon\otimes 1))\\
&=2,
\end{align*}
hence $[f,f]\neq 0$.
\par
On the other hand, if $A=B\oplus I$ as associative $\K$-algebra,
then the $L_{\infty}$ structure on $K$ is trivial. Indeed,
considering $K[1]$ as a DGLA with trivial bracket, the obvious map
\[ K[1]=\Hoch^\bullet(I,B)\to \Hoch^\bullet(A,A)\]
gives a commutative diagram of morphisms of DGLAs
\[\begin{matrix}
0&\to&\ker\alpha\\
\mapver{}&&\mapver{\chi}\\
K[1]&\to&\Hoch^\bullet(A,A)
\end{matrix}\]
such that the composition $K\to C_{\chi}\to K$ is the identity.
Therefore the $L_\infty$-algebra structure induced on $K$ is
isomorphic to $\widetilde{C}(0\hookrightarrow K[1])$, which is a
trivial $L_\infty$-algebra.
\end{example}

\end{document}